\documentclass{amsart}
\begin{document}
\allowdisplaybreaks
\parskip10pt
\parindent0pt
\def\e#1{\mathbf{e}_{#1}}
\def\RR{\mathbb{R}^2}
\def\R{\mathbb{R}}
\def\H{H}
\def\<{\left<}
\def\>{\right>}
\def\a{\mathbf{a}}
\def\b{\mathbf{b}}
\def\v{\mathbf{v}}
\def\w{\mathbf{w}}
\def\x{\mathbf{x}}
\def\ll#1{\lambda_{#1}}
\def\mon{monochromatic}
\def\hom{homothetic}
\def\ie{{\it i.e.}}
\def\bfd#1{{\bf#1}}
\def\A(#1,#2){$\underline{\phi}(#1,#2)$}
\def\B(#1,#2,#3){$\underline{\delta}(#1,#2,#3)$}
\newtheorem{theorem}{{\bf Theorem}}
\title{Gallai's Theorem}
\author{R.~D.~Maddux}
\address{Dept. Math., Iowa State University, Ames, Iowa 50011}
\date{Oct 14, 2013}
\begin{abstract}
This paper presents a proof of Gallai's Theorem, adapted from A.~Soifer's presentation
in~\emph{The Mathematical Coloring Book} \cite{Soifer} of E.~Witt's 1952 proof of
Gallai's Theorem \cite{Witt}.
\end{abstract}
\maketitle
\section{Introduction}
Gallai's Theorem states that if the points in the Euclidean plane are colored with
finitely many colors, then for every finite subset of the plane there is a monochromatic
homothetic copy of that set. A ``homothetic copy'' of a set is its image under first a
dilation and then a translation, while ``monochromatic'' means that every point in the
copy receives the same color. Another way to state Gallai's Theorem is that for every
finite subset $X$ of the plane and every number $k$ of colors, there is another finite
subset $Y$ of the plane, such that if $Y$ is colored by $k$ colors then $Y$ contains a
\mon\ \hom\ image of $X$. (The two versions are equivalent by the compactness theorem
for first-order logic.)  I first heard about Gallai's Theorem in late 2012 from Jeremy
Alm, who wrote and sent me a paper about extending Gallai's Theorem. I became intrigued
by the theorem, tried to prove it, and eventually consulted the source referenced in
Jeremy's paper, Chapter~42 of A.~Soifer's~\emph{The Mathematical Coloring
Book}~\cite{Soifer}.  Soifer~\cite{Soifer} presents an expanded account of a proof of
Gallai's Theorem published by Ernst Witt in 1952 \cite{Witt}. Soifer thought Witt's
proof was incomprehensibly brief, so he added details. Soifer's proof seemed
incomprehensibly brief to me, so I added even more details.
\section{Definitions}
Consider an arbitrary but fixed infinite sequence
$S=\<\e0,\e1,\e2,\cdots,\e{i},\cdots\>$ of distinct points in $\RR$, where
$\e0=(0,0)$. For every integer $n\geq1$ let
\begin{equation*}
S_n:=\{\e0,\cdots,\e{n-1}\}.
\end{equation*}
For example, $S_1=\{\e0\}$, $S_2=\{\e0,\e1\}$, and $S_3=\{\e0,\e1,\e2\}$.  For every set
$X$, $|X|$ is the number of elements of $X$, called the \bfd{cardinality} of $X$. The
cardinality of $S_n$ is $n$, \ie, $|S_n|=n$. For every positive integer $k$,
a~\bfd{$k$-element set} is a set $X$ such that $|X|=k$.  A function from a set $Y$ to a
$k$-element set $X$ is a \bfd{$k$-coloring} of $Y$.

Let $\H$ be the set of all functions $h:\RR\to\RR$ such that, for some \bfd{displacement
vector} $\a\in\RR$ and some positive \bfd{dilation scalar} $\lambda\in\R$, we have
$h(\v)=\a+\lambda \v$ for every $\v\in\RR$, where addition and scalar multiplication are
performed component-wise.  $H$ is the set of \bfd{homotheties} from $\RR$ to $\RR$,
consisting of all compositions of translations and dilations, \ie,
\begin{equation*}
	H:=\{h:(\exists\a\in\RR)(\exists\lambda\in\R)
	(0\leq\lambda\land(\forall\v\in\RR)(h(\v)=\a+\lambda \v))\}.
\end{equation*}
If $V\subseteq\RR$ and $h\in H$, then $h(V)$ is the \bfd{image of} $V$ \bfd{under} $h$,
\ie,
\begin{equation*}
	h(V):=\{h(\v):\v\in\RR\}.
\end{equation*}
For any set $V\subseteq\RR$ and any $n\geq2$, $E_n(V)$ is the set consisting of the
images of $S_n$ under those homotheties that map $S_{n-1}$ into $V$, \ie,
\begin{equation*}
	E_n(V):=\bigcup\{h(S_n):h\in\H,h(S_{n-1})\subseteq V\}.
\end{equation*}
For any vector $\v\in\RR$ and any two sets $V,W\subseteq\RR$, the \bfd{complex sum} of
$\v$, or $V$, and $W$ is the set $\v+W$, or $V+W$, of vectors obtained by adding $\v$,
or any vector from $V$, to any vector from $W$, that is,
\begin{align*}
	\v+W&:=\{\v+\w:\w\in W\},
\\	V+W&:=\{\v+\w:\v\in V,\w\in W\}=\bigcup_{\v\in V}(\v+W).
\end{align*}
For all integers $k\geq2$, $n\geq3$, and $m\geq1$, define sets
$\Phi(n,k),\Delta(n,k,m)\subseteq\RR$ as follows.
\begin{align}
  \label{1}	\Phi(2,k)	&:=\{\e0,\e1,2\e1,3\e1,\cdots,(k-1)\e1,k\e1\}
\\\label{2}	\Delta(n,k,1)	&:=E_n(\Phi({n-1},k))
\\\label{3}	\Delta(n,k,m+1)	&:=\Delta(n,k^{|\Delta(n,k,m)|},1)+\Delta(n,k,m)
\\\label{4}	\Phi(n,k)	&:=\Delta(n,k,k)
\end{align}
For example, with $n=3$, $k=2$, and $m=1$, we have
\begin{align*}
  \Phi(2,2)
  &=\{\e0,\e1,2\e1\}
  &&\text{by \eqref{1}}
\intertext{There are three homotheties $h_1,h_2,h_3\in H$ that map $S_2$ into
$\{\e0,\e1,2\e1\}$, defined for all $\v\in\RR$ by $h_1(\v)=\v,$ $h_2(\v)=2\v,$ and
$h_3(\v)=\e1+\v.$ Hence $|E_3(\{\e0,\e1,2\e1\})|\leq6$ since}
	E_3(\{\e0,\e1,2\e1\})
  &=	\bigcup\{h(S_3):h\in\H,h(S_2)\subseteq\{\e0,\e1,2\e1\}\}
\\&=	\bigcup\{h(S_3):h\in\{h_1,h_2,h_3\}\}
\\&=	h_1(S_3)\cup h_2(S_3)\cup h_3(S_3)
\\&=	S_3\cup h_2(\{\e0,\e1,\e2\})\cup h_3(\{\e0,\e1,\e2\})
\\&=	\{\e0,\e1,\e2\}\cup\{2\e0,2\e1,2\e2\}\cup\{\e0+\e1,2\e1,\e1+\e2\}
\\&=	\{\e0,\e1,\e2,2\e1,2\e2,\e1+\e2\}.
\intertext{If $|E_3(\{\e0,\e1,2\e1\})|=6$ then}
	\Phi(3,2)
  &=	\Delta(3,2,2)
	&&\text{by \eqref{4}}
\\&=	\Delta(3,2^{|\Delta(3,2,1)|},1)+\Delta(3,2,1)
	&&\text{by \eqref{3}}
\\&=	\Delta(3,2^{|E_3(\Phi(2,2))|},1)+E_3(\Phi(2,2))
	&&\text{by \eqref{2}}
\\&=	\Delta(3,2^{|E_3(\{\e0,\e1,2\e1\})|},1)+E_3(\{\e0,\e1,2\e1\})
	&&\text{by \eqref{1}}
\\&=	\Delta(3,2^6,1)+\{\e0,\e1,2\e1,\e2,2\e2,\e1+\e2\}	
\\&=	E_3(\Phi(2,64))+\{\e0,\e1,2\e1,\e2,2\e2,\e1+\e2\}
	&&\text{by \eqref{3}}
\end{align*}
Since $\Phi(2,64)$ has 65 vectors by \eqref{1}, there are
$\begin{pmatrix}64\\2\end{pmatrix}=2080$ homotheties mapping $S_2$ into $\Phi(2,64)$,
hence $E_3(\Phi(2,64))$ may contain as many as 65+2080=2145 vectors, and $\Phi(3,2)$ may
contain as many as 2145+6=2151.

For integers $n,k,m$, let \A(n,k) be the statement that 
\begin{quote}
for every $k$-coloring of $\Phi(n,k)$, $\Phi(n,k)$ contains a monochromatic homothetic
image of $S_n$,
\end{quote}
and let \B(n,k,m) be the statement that 
\begin{quote}
for every $k$-coloring $f$ of $\Delta(n,k,m)$ there are a vector $\a\in\RR$ and scalars
$\ll0,\ll1,\cdots,\ll{m}\in\R$ such that $0=\ll0<\ll1<\cdots<\ll{m}$ and if
$h_{i,j}(\v)=\a+\ll{i}\e{n}+(\ll{j}-\ll{i})\v$ whenever $\v\in\RR$ and $0\leq i<j\leq
m$, then $h_{i,j}(S_n)\subseteq\Delta(n,k,m)$ and $|f(h_{i,j}(S_{n-1}))|=1$.
\end{quote}
\begin{theorem}[Gallai's Theorem]
If $3\leq n\in\omega$, $2\leq k\in\omega$, and $1\leq m\in\omega$, then~\A(n,k)
and~\B(n,k,m).
\end{theorem}
The proof is by induction. For each of \eqref{1}--\eqref{4}, we show that if the sets on
the right side have the requisite properties, then the set on left side does so as
well. here is the logic of the proof:
\bgroup
\def\A(#1,#2){\underline{\phi}(#1,#2)}
\def\B(#1,#2,#3){\underline{\delta}(#1,#2,#3)}
\begin{align*}						&\A(2,k)&&\text{Base Cases}
\\\A(n-1,k)\implies					&\B(n,k,1)		
\\\B(n,k,m)\land\B(n,k^{|\Delta(n,k,m)|},1)\implies	&\B(n,k,m+1)		
\\\B(n,k,k)\implies					&\A(n,k)		
\end{align*}
\egroup
\subsection*{Proof of \A(2,k)}
Since $\e0=(0,0)$, \eqref{1} implies that
$$	\Phi(2,k)=\{0\e1,\e1,2\e1,\cdots,k\e1\}.
$$
We must show that for every $k$-coloring $f$ of $\Phi(2,k)$,
$\{0\e1,\e1,2\e1,\cdots,k\e1\}$ contains a \mon\ \hom\ image of $S_2=\{\e0,\e1\}$.  The
number of vectors in $\Phi(2,k)$ is $k+1$, one more than the number of colors, so there
are (at least) two vectors in $\Phi(2,k)$ that get the same color, say $f(i\e1)=f(j\e1)$
for some $i,j$ such that $0\leq i<j\leq k$. Define $h\in H$ by $h(\v)=i\e1+(j-i)\v$ for
all $v\in\RR$. Then $h$ maps $S_2=\{\e0,\e1\}$ to $\{i\e1,j\e1\}\subseteq\Phi(2,k)$
because
\begin{align*}
	h(\e0)&=i\e1+(j-i)\e0=i\e1+(j-i)(0,0)=i\e1,
\\	h(\e1)&=i\e1+(j-i)\e1=i\e1+j\e1-i\e1=j\e1.
\end{align*}
This \hom\ image is monochromatic by the assumption $f(i\e1)=f(j\e1)$.
\subsection*{Proof of \B(n,k,1) from $n\geq3$ and \A(n-1,k)}
Consider a $k$-coloring $f$ of $\Delta(n,k,1)$.  By the inductive assumption \A(n-1,k),
$\Phi({n-1},k)$ contains a \mon\ \hom\ image of $S_{n-1}$, so there is some $h\in H$
such that
\begin{align}	\label{e}
	h(S_{n-1})\subseteq\Phi({n-1},k),
\\		\label{f}
	|f(h(S_{n-1}))|=1.
\end{align}
From \eqref{e}, the definition of $E_n$, and definition \eqref{2}, it follows that
\begin{equation}\label{g}
	h(S_n)\subseteq E_n(\Phi({n-1},k))=\Delta(n,k,1).
\end{equation}
Now, according to \B(n,k,1), we must find $\a\in\RR$ and $\ll0,\ll1\in\R$ such that
$0=\ll0<\ll1$, $h_{0,1}(S_n)\subseteq\Delta(n,k,1)$, and $|f(h_{0,1}(S_{n-1}))|=1$,
where $$h_{0,1}(\v)=\a+\ll0\e{n}+(\ll1-\ll0)\v=\a+\ll1\v$$ whenever $\v\in\RR$.  It
suffices to let $\a$ and $\ll1$ be the displacement vector and dilation scalar
associated with $h$, for then $h_{0,1}=h$ and the two desired equations are \eqref{f}
and \eqref{g}.
\begin{figure}\setlength{\unitlength}{0.5mm}\boxed{
\begin{picture}(263,306)
\put(0,0){\circle*{2.1}$\e0$}
\put(10,0){\circle*{2.1}$\e1$}
\put(10,5){\circle*{2.1}$\e2$}
\put(0,13){\circle*{2.1}$\e3$}
\put(13,120){\circle*{2.1}}
\put(13,146){\circle*{2.1}}
\put(13,159){\circle*{2.1}}
\put(13,185){\circle*{2.1}}
\put(23,146){\circle*{2.1}}
\put(23,151){\circle*{2.1}}
\put(33,120){\circle*{2.1}}
\put(33,130){\circle*{2.1}}
\put(33,159){\circle*{2.1}}
\put(33,169){\circle*{2.1}}
\put(43,120){\circle*{2.1}}
\put(43,135){\circle*{2.1}}
\put(43,146){\circle*{2.1}}
\put(43,161){\circle*{2.1}}
\put(63,120){\circle*{2.1}}
\put(63,145){\circle*{2.1}}
\put(123,140){\circle{2.1}}
\put(123,166){\circle{2.1}}
\put(123,179){\circle{2.1}}
\put(123,205){\circle{2.1}}
\put(123,296){\circle{2.1}}
\put(133,166){\circle{2.1}}
\put(133,171){\circle{2.1}}
\put(143,140){\circle{2.1}}
\put(143,150){\circle{2.1}}
\put(143,179){\circle{2.1}}
\put(143,189){\circle{2.1}}
\put(153,140){\circle{2.1}}
\put(153,155){\circle{2.1}}
\put(153,166){\circle{2.1}}
\put(153,181){\circle{2.1}}
\put(173,140){\circle{2.1}}
\put(173,165){\circle{2.1}}
\put(193,205){\circle{2.1}}
\put(193,240){\circle{2.1}}
\put(213,179){\circle{2.1}}
\put(213,224){\circle{2.1}}
\put(223,166){\circle{2.1}}
\put(223,216){\circle{2.1}}
\put(243,140){\circle{2.1}}
\put(243,200){\circle{2.1}}
\put(110,20){\circle{2.1}}
\put(180,20){\circle{2.1}}
\put(180,55){\circle{2.1}}
\put(13,120){\circle*{2.1}$\a$ }
\put(63,120){\circle*{2.1}$\a+\ll3\e1$}
\put(63,145){\circle*{2.1}$\a+\ll3\e2$}
\put(13,185){\circle*{2.1}$\a+\ll3\e3$}
\put(123,205){\circle{2.1}$\b+\a+\ll3\e3$}
\put(193,205){\circle{2.1}$\b+\mu_1\e1+\a+\ll3\e3$}
\put(193,240){\circle{2.1}$\b+\mu_1\e2+\a+\ll3\e3$}
\put(123,296){\circle{2.1}$\b+\mu_1\e3+\a+\ll3\e3$}
\put(110,20){\circle*{2.1}$\b$}
\put(180,20){\circle*{2.1}$\b+\mu_1\e1$}
\put(180,55){\circle*{2.1}$\b+\mu_1\e2$}
\put(110,111){\circle*{2.1}$\b+\mu_1\e3$}
\end{picture}}
\label{fig1}
\caption{$m=n=4$, $\e0=(0,0)$, $\e1=(10,0)$, $\e2=(10,5)$, $\e3=(0,13)$, $\a=(13,120)$,
$\b=(110,20)$, $\mu_1=7$, $\ll1=2$, $\ll2=3$, $\ll3=5$, $\ll4=5+\mu_1$.}
\end{figure}
\subsection*{Proof of \B(n,k,m+1) from \B(n,k,m) and \B(n,k^{|\Delta(n,k,m)|},1)}
Suppose $f$ is a $k$-coloring of $\Delta(n,k,m+1)$. By definition \eqref{3}, $f$ assigns
a color to every vector obtained by adding a vector from
$\Delta(n,k^{|\Delta(n,k,m)|},1)$ to a vector in $\Delta(n,k,m)$. Therefore, for every
$\v\in\Delta(n,k^{|\Delta(n,k,m)|},1)$, we may let $f'(\v)$ be the $k$-coloring of
$\Delta(n,k,m)$ defined by
\begin{equation}\label{5}
	f'(\v)(\w)=f(\v+\w)\quad\text{for every }\w\in\Delta(n,k,m).
\end{equation}
This gives us a new coloring $f'$ that assigns each vector
$\v\in\Delta(n,k^{|\Delta(n,k,m)|},1)$ to an element $f'(\v)$ of a set whose cardinality
is $k^{|\Delta(n,k,m)|}$.  Thus $f'$ is a $k^{|\Delta(n,k,m)|}$-coloring of
$\Delta(n,k^{|\Delta(n,k,m)|},1)$.  From the inductive hypothesis
\B(n,{k^{|\Delta(n,k,m)|}},1), applied to the coloring $f'$, we know there are
$\b\in\RR$ and $\mu_0,\mu_1\in\R$ such that $0=\mu_0<\mu_1$ and, assuming
\begin{equation}\label{6}
	g_{0,1}(\v)=\b+\mu_0\e{n}+(\mu_1-\mu_0)\v=\b+\mu_1\v
	\quad\text{for all $\v\in\RR$,}
\end{equation}
we have
\begin{align}
  &\label{7}
	g_{0,1}(S_n)\subseteq\Delta(n,k^{|\Delta(n,k,m)|},1),
\\&\label{7a}
	|f'(g_{0,1}(S_{n-1}))|=1.
\end{align}
From the inductive hypothesis \B(n,k,m), applied to the $k$-coloring $f'(\b)$, we know
there are $\a\in\RR$ and $\ll0,\ll1,\cdots,\ll{m}\in\R$ such that
$0=\ll0<\ll1<\cdots<\ll{m}$ and, assuming
\begin{equation}\label{9}
	h_{i,j}(\v)=\a+\ll{i}\e{n}+(\ll{j}-\ll{i})\v
	\quad\text{if $\v\in\RR$ and $0\leq i<j\leq m$,}
\end{equation}
we have
\begin{align}	&\label{10}
	h_{i,j}(S_n)\subseteq\Delta(n,k,m),
\\&\label{11}
	|f'(\b)(h_{i,j}(S_{n-1}))|=1.
\end{align}
This situation is illustrated in Fig.~\ref{fig1} in case $m=n=1$. Next we prove that if
$\w\in\Delta(n,k,m)$, $0\leq l<n-1$, and $0\leq i<m$, then
\begin{align}
  \label{8}
	f(\b+\w)
&=	f(\b+\mu_1\e{l}+\w),
\\\label{12}
	f(\b+\a+\ll{i}\e{n})
&=	f(\b+\a+\ll{i}\e{n}+(\ll{m}-\ll{i})\e{l}).
\end{align}
Proof of \eqref{8}: If $0\leq l<n-1$, then $\e0,\e{l}\in S_{n-1}$ and, by \eqref{6},
$g_{0,1}(\e0)=\b$ and $g_{0,1}(\e{l})=\b+\mu_1\e{l}$, hence $\b,\b+\mu_1\e{l}\in
g_{0,1}(S_{n-1})$ so from~\eqref{7a} we conclude that $f'(\b)=f'(\b+\mu_1\e{l})$. This
holds whenever $0\leq l<n-1$, so by~\eqref{5}, we obtain \eqref{12}.

\par\noindent Proof of \eqref{12}: If $0\leq l<n-1$ and $0\leq i<m$ then $\e0,\e{l}\in
S_{n-1}$ and, by~\eqref{9} (applied with $0\leq i<j=m$), we have
\begin{align*}
	\a+\ll{i}\e{n}&=h_{i,m}(\e0)\in h_{i,m}(S_{n-1}),
\\	\a+\ll{i}\e{n}+(\ll{m}-\ll{i})\e{l}&=h_{i,m}(\e{l})\in h_{i,m}(S_{n-1}).
\end{align*}
These two vectors are in $\Delta(n,k,m)$ by \eqref{10} and get the same color from
$f'(\b)$ by \eqref{11} with $j=m$, \ie,
$$	f'(\b)(	\a+\ll{i}\e{n})=f'(\b)(	\a+\ll{i}\e{n}+(\ll{m}-\ll{i})\e{l}).
$$
The conclusion follows by~\eqref{5} with $\v=\b$.

Define a new system $h'$ as follows.  Let $\ll{m+1}=\mu_1+\ll{m}$ and let
\begin{equation}\label{13}
	h'_{i,j}(\v)=\b+\a+\ll{i}\e{n}+(\ll{j}-\ll{i})\v
	\quad\text{if $\v\in\RR$ and $0\leq i<j\leq m+1$.}
\end{equation}
We must now show for this new system that $0=\ll0<\ll1<\cdots<\ll{m}<\ll{m+1}$, (which
we get immediately from $0<\mu_1$ and $\ll{m+1}=\mu_1+\ll{m}>\ll{m}$) and, whenever
$\v\in\RR$ and $0\leq i<j\leq m+1$, we have
\begin{align}
  \label{14}&	h'_{i,j}(S_n)\subseteq\Delta(n,k,m+1),
\\\label{15}&	|f(h'_{i,j}(S_{n-1}))|=1.
\end{align}
\par\noindent{\sc CASE~1.} $0\leq i<j\leq m$.
\begin{align*}
	h'_{i,j}(S_n)
	&=\b+h_{i,j}(S_n)				&&\text{by \eqref{9}, \eqref{13}}
\\	&\subseteq\b+\Delta(n,k,m)			&&\text{by \eqref{10}}
\\	&=g_{0,1}(\e0)+\Delta(n,k,m)			&&\text{by \eqref{6}, $\e0=(0,0)$}
\\	&\subseteq g_{0,1}(S_n)+\Delta(n,k,m)		&&\e0\in S_n
\\	&\subseteq\Delta(n,k^{|\Delta(n,k,m)|},1)
		+\Delta(n,k,m)				&&\text{by \eqref{7}}
\\	&=\Delta(n,k,m+1)				&&\text{by \eqref{3}}
\end{align*}
so \eqref{14} holds in this case, and
\begin{align*}
	|f(h'_{i,j}(S_{n-1}))|
	&=|f(\b+h_{i,j}(S_{n-1}))|	&&\text{by \eqref{9}, \eqref{13}}
\\	&=|f'(\b)(h_{i,j}(S_{n-1}))|	&&\text{by \eqref{5}}
\\	&=1				&&\text{by \eqref{11}}
\end{align*}
so \eqref{15} also holds in this case.
\par\noindent{\sc CASE~2.} $0\leq i<j=m+1$.
Proof of \eqref{14}: If $\v\in S_n$, then
\begin{align*}
	h'_{i,m+1}(\v)
  &=	\b+\a+\ll{i}\e{n}+(\ll{m+1}-\ll{i})\v		&&\text{\eqref{13}}
\\&=	\b+\a+\ll{i}\e{n}+(\ll{m}+\mu_1-\ll{i})\v
\\&=	\b+\mu_1\v+\a+\ll{i}\e{n}+(\ll{m}-\ll{i})\v
\\&=	g_{0,1}(\v)+h_{i,m}(\v)				&&\text{\eqref{6}, \eqref{9}}
\\&\in	\Delta(n,k^{|\Delta(n,k,m)|},1)+\Delta(n,k,m)	&&\text{\eqref{7}, \eqref{10}}
\\&=	\Delta(n,k,m+1)					&&\text{\eqref{3}}
\end{align*}
For \eqref{15}, consider an arbitrary $l$ with $0\leq l<n-1$, \ie, any arbitrary
$\e{l}\in S_{n-1}$.  Note that by~\eqref{9} and~\eqref{10} with $j=m$ we have
\begin{equation}\label{p}
	\a+\ll{i}\e{n}+(\ll{m}-\ll{i})\e{l}=h_{i,m}(\e{l})\in\Delta(n,k,m).
\end{equation}
Hence
\begin{align*}
	f(h'_{i,m+1}(\e0))
	&=f(\b+\a+\ll{i}\e{n})
	&&\text{by \eqref{13}}
\\	&=f(\b+\a+\ll{i}\e{n}+(\ll{m}-\ll{i})\e{l})
	&&\text{by \eqref{12} if $i<m$, triv if $i=m$}
\\	&=f(\b+\mu_1\e{l}+\a+\ll{i}\e{n}+(\ll{m}-\ll{i})\e{l})
	&&\text{by \eqref{p}, \eqref{8}}
\\	&=f(\b+\a+\ll{i}\e{n}+(\mu_1+\ll{m}-\ll{i})\e{l})
	&&\text{by computation}
\\	&=f(\b+\a+\ll{i}\e{n}+(\ll{m+1}-\ll{i})\e{l})
	&&\text{by definition of $\ll{m+1}$}
\\	&=f(h'_{i,m+1}(\e{l}))
	&&\text{by \eqref{13}}
\end{align*}
This argument shows that \eqref{15} holds in this case.
\subsection*{Proof of  \A(n,k) from \B(n,k,k)}
Since $\Phi(n,k)=\Delta(n,k,k)$, we need only show that $\Delta(n,k,k)$ is a set which,
for every $k$-coloring of $\Delta(n,k,k)$, contains a \mon\ \hom\ image of $S_n$.
Suppose $f$ is a $k$-coloring of $\Delta(n,k,k)$.  From the inductive hypothesis
\B(n,k,k), we know there are $\a\in\RR$ and $\ll0,\ll1,\cdots,\ll{k}\in\R$ such
that $0=\ll0<\ll1<\cdots<\ll{k}$ and if
\begin{equation*}
	h_{i,j}(\v)=\a+\ll{i}\e{n}+(\ll{j}-\ll{i})\v
\end{equation*}
whenever $\v\in\RR$ and $0\leq i<j\leq k$, then $h_{i,j}(S_n)\subseteq\Delta(n,k,k)$ and
$|f(h_{i,j}(S_{n-1}))|=1$.  Let's look at the images of $\e{n}$ and $\e0$ under all
these homotheties $h_{i,j}$ with $0\leq i<j\leq k$. We have
\begin{equation*}
	h_{i,j}(\e{n})=\a+\ll{i}\e{n}+(\ll{j}-\ll{i})\e{n}=\a+\ll{j}\e{n},
\end{equation*}
and, since $\e0=(0,0)$, 
\begin{equation*}
	h_{i,j}(\e0)=\a+\ll{i}\e{n}+(\ll{j}-\ll{i})\e0=\a+\ll{i}\e{n}.
\end{equation*}
The images of $\e{n}$ are $\a+\ll1\e{n},\cdots,\a+\ll{k-1}\e{n},\a+\ll{k}\e{n}$, and,
since $\ll0=0$, the images of $\e0$ are $\a,\a+\ll1\e{n},\cdots,\a+\ll{k-1}\e{n}$.  The
union of these two lists, namely
$$	\{\a,\a+\ll1\e{n},\cdots,\a+\ll{k-1}\e{n},\a+\ll{k}\e{n}\}
$$
is a $(k+1)$-element set, to which are assigned only $k$ colors. Two elements receive
the same color, hence there are $i,j$ such that $0\leq i<j\leq k$ and
\begin{equation*}
	f(\a+\ll{i}\e{n})=f(\a+\ll{j}\e{n}).
\end{equation*}
However, $\a+\ll{i}\e{n}=h_{i,j}(\e0)$ and $\a+\ll{j}\e{n}=h_{i,j}(\e{n})$, so
\begin{equation*}
	f(h_{i,j}(\e0))=f(h_{i,j}(\e{n})).
\end{equation*}
Of course, we already know $|f(h_{i,j}(S_{n-1}))|=1$, so this last equation tells us,
since $\e0\in S_{n-1}$, that $\e{n}$ also receives the same color as all the other
elements of $S_{n-1}$. But $S_n=S_{n-1}\cup\{\e{n}\}$, so in fact we have shown
$|f(h_{i,j}(S_n))|=1$. We also know that $h_{i,j}(S_n)\subseteq\Delta(n,k,k)$, so we
conclude that $\Delta(n,k,k)$ does indeed contain a \mon\ \hom\ image of $S_n$.

\end{document}